\documentclass[10pt,legno]{article}
\usepackage{amsmath, amssymb}
\usepackage{latexsym}
\usepackage{enumerate}
\begin{document}
\newcommand{\qed}{\hfill \ensuremath{\square}}
\newtheorem{thm}{Theorem}[section]
\newtheorem{cor}[thm]{Corollary}
\newtheorem{lem}[thm]{Lemma}
\newtheorem{prop}[thm]{Proposition}
\newtheorem{defn}[thm]{Definition}
\newcommand{\proof}{\vspace{1ex}\noindent{\em Proof}. \ }
\newtheorem{pro}[thm]{proof}
\newtheorem{ide}[thm]{Idee}
\newtheorem{rem}[thm]{Remark}
\newtheorem{ex}[thm]{Example}
\bibliographystyle{plain}
\numberwithin{equation}{section}
\numberwithin{equation}{section}
\newcounter{saveeqn}
\newcommand{\subeqn}{\setcounter{saveeqn}{\value{equation}}%
 \stepcounter{saveeqn}\setcounter{equation}{0}
\renewcommand{\theequation}{\mbox{\arabic{section}.\arabic{saveeqn}\alph{=
equation}}}} 

\newcommand{\reseteqn}{\setcounter{equation}{\value{saveeqn}}%
\renewcommand{\theequation}{\arabic{section}.\arabic{equation}}}

\def\nm{\noalign{\medskip}}
\newcommand{\Om}{\Omega}
\newcommand{\om}{\omega}
\newcommand{\Real}{\mathbb{R}}
\newcommand{\nuu}{\tilde{\nu}}
\newcommand{\bohm}{{\partial}{\ohm}}
\newcommand{\la}{\langle}
\newcommand{\ra}{\rangle}
\newcommand{\ms}{\mathcal{S}_\ohm}
\newcommand{\mk}{\mathcal{K}_\ohm}
\newcommand{\mks}{\mathcal{K}_\ohm ^{\ast}}
\newcommand{\grad}{\bigtriangledown}
\newcommand{\ds}{\displaystyle}
\newcommand{\pf}{\medskip \noindent {\sl Proof}. ~ }
\newcommand{\p}{\partial}
\renewcommand{\a}{\alpha}
\newcommand{\z}{\zeta}
\newcommand\q{\quad}
\newcommand{\pd}[2]{\frac {\p #1}{\p #2}}
\newcommand{\pdl}[2]{\frac {\p^2 #1}{\p #2}}
\newcommand{\dbar}{\overline \p}
\newcommand{\eqnref}[1]{(\ref {#1})}
\newcommand{\na}{\nabla}
\newcommand{\ep}{\epsilon}
\newcommand{\vp}{\varphi}
\newcommand{\fo}{\forall}
\newcommand{\Tcal}{\mathcal{T}}
\newcommand{\ta}{\theta}
\newcommand{\Scal}{\mathcal{S}}
\newcommand{\Dcal}{\mathcal{D}}
\newcommand{\Kcal}{\mathcal{K}}
\newcommand{\Acal}{\mathcal{A}}
\newcommand{\Ecal}{\mathcal{E}}
\newcommand{\Ncal}{\mathcal{N}}
\newcommand{\Abar}{\overline A}
\newcommand{\Rcal}{\mathcal{R}}
\newcommand{\Lcal}{\mathcal {L}}
\newcommand{\Gcal}{\mathcal {G}}
\newcommand{\Cbar}{\overline C}
\newcommand{\Ebar}{\overline E}
\newcommand{\RR}{\mathbb{R}}
\newcommand{\CC}{\mathbb{C}}
\newcommand{\NN}{\mathbb{N}}
\newcommand{\Z}{\mathbb{Z}}

\title{ Sensitivity Analysis of  the  Current, the Electrostatic  Capacity, and the Far-Field of the Potential     with Respect to Small Perturbations in the Surface  of   a Conductor}
\date{}
\author{ Jihene Lagha \thanks{ Universit\'e de Tunis El Manar, Facult\'e des Sciences de Tunis, LR11ES13 Laboratoire d'Analyse stochastique et Applications, 2092 Tunis,  Tunisie (lagha.jihene@yahoo.fr)}
\and Habib Zribi\thanks{Department of Mathematics, College of Science,
University of Hafr Al Batin, P.O. 1803, Hafr Al Batin 31991, Saudi Arabia (zribi.habib@yahoo.fr)}}
\maketitle

\begin{abstract}
We derive  asymptotic  expansions of  the  current,  the  electrostatic  capacity, and  the far-field of the  electrostatic potential resulting  from   small perturbations in the  shape  of an isolated   conductor with  $\mathcal{C}^2$-surface.
Our derivation is rigorous by using systematic way, based  on layer potential techniques and the field expansion (FE) method.  We then use these results  to study    the sensibility analysis of the first  eigenvector  of the $L^2$-adjoint   of the  Neumann-Poincar\'e  (NP)   operator with respect to small perturbations in the surface  of   its domain.
\end{abstract}

\noindent {\footnotesize Mathematics Subject Classification
(MSC2000): 35B30,   35B40}

\noindent {\footnotesize Keywords:  Isolated conductor, electrostatic  capacity,
small surface perturbations,  boundary integral method,
field expansion method,  Laplace equation, Neumann-Poincar\'e operator type

}

\section{Introduction and statement of the  main results}
Suppose that an isolated   conductor   occupies a bounded domain $\Om$ in $\RR^3,$ with a connected $\mathcal{C}^2$-surface $\p \Om$. A conductor is a volume which contains free charges. In the presence of an external electric field,
each free charge in the conductor redistributes and very quickly reaches electrostatic equilibrium  (see  \cite{volume2, Maxwell}).
 The free charges are redistributed in such a way the electric field  inside the conductor vanishes and the
 electrical filed is always perpendicular everywhere on the surface of the conductor. The electrostatic potential is constant throughout
 the volume of the conductor and has the same value on its surface, let's say $1$ volt.  More precisely,  let     $u$  be the   electrostatic potential   in the presence of a conductor $\Om$  at equilibrium in $\RR^3$. It is  the unique solution of the following problem
\begin{equation}\label{u}
\left\{
  \begin{array}{ll}
\ds  \Delta u=0 \quad \mbox{ in } \RR^3 \backslash \overline{\Om}, \\
\nm  \ds u=1\q \quad \mbox{on  }  \p \Om,\\
\nm \ds  \lim_{|x|\rightarrow \infty}u(x) =0.
  \end{array}
\right.
\end{equation}
The electrostatic capacity  with respect to  infinity of  the conductor $\Om$, denoted $\textrm{cap}(\Om)$,  is
defined as the ratio of the charge in  equilibrium on it to the value of the potential $u$ at its  surface $\p \Om$. That is,   the capacity
is the charge producing this potential which is given by Gauss' integral (see for instance \cite{Poincare, Polya-Szego-book, Kellogg})
\begin{align}\label{capacity-Om}
\textrm{cap}(\Om)
=-\frac{1}{4\pi}\int_{\p\Om} \frac{\p u}{\p n}(x)  d\sigma(x)\q \Big(=-\frac{1}{4\pi}\int_{\p\Om} \frac{\p u}{\p n}(x) u(x)  d\sigma(x)\Big),
\end{align}
where  $n$  and $d\sigma$  are the unit outward normal  and the length element to the boundary $\p \Om$, respectively.
The electrostatic capacity  $\textrm{cap}(\Om)$  may also  be defined as the quantity of  electrical charge
which  must be given to the conductor $\Om$ to raise its potential to the value unity, it  depends on its own form and size; being greater as the seize increased. In this work, we investigate the sensitivity analysis of the electrostatic capacity  with small changes in the form of its domain.

It is well-known that  $ 0 \leq \textrm{cap}(\Om)<
+\infty$, which can be easily  proved
  by applying  the Green's identity to the  integral in \eqref{capacity-Om} over the unbounded domain $\RR^3\backslash \overline{\Om}$.  The electrostatic capacity can be determined from the far-field of the  electrostatic potential  $u$  defined in  \eqref{u}. In fact, in \cite{ Polya, Jerison}, $u$  has the asymptotic expansion
\begin{align}\label{asymptotic-infinity}
\ds u(x)=\frac{\textrm{cap}(\Om)}{ |x|} +O(\frac{1}{|x|^{2}})\quad \mbox{as } |x|\rightarrow \infty.
\end{align}
Therefore, in order to find  the electrostatic capacity, we have to pick out the coefficient of $1/|x|$ in the
 asymptotic expansion \eqref{asymptotic-infinity}.  The capacities  are known analytically for a few simple shapes like   sphere, ellipsoid, lens, spindle, and  anchor-ring. See  \cite{Szego,Polya-Szego-book}.

Let  $\Om_\ep$ be an $\ep$-perturbation of $\Om$,  {\it i.e.}, there is a
function $h\in C^{1}(\p \Om)$ such that $\p \Om_\ep$ is given by
\begin{align*}
\ds \partial {\Om_\ep}=\big \{\tilde{x}=x+\ep h(x) n(x):=\Psi_\ep (x) | x \in \p \Om
\big \}.
\end{align*}
We denote  by  $u_{\ep}$   the perturbed  electrostatic potential  in the presence of the conductor  $\Om_\ep$ in electrostatic equilibrium.  It is  the unique  solution
 of the following problem
\begin{equation}\label{uep}
\left\{
  \begin{array}{ll}
\ds \Delta u_{\ep}=0 \quad \mbox{ in } \RR^3 \backslash \overline{\Om}_\ep, \\
\nm \ds u_{\ep}=1 \q\quad \mbox{on  }  \p \Om_\ep,\\
\nm \ds  \lim_{|x|\rightarrow \infty}u_{\ep}(x) =0.
  \end{array}
\right.
\end{equation}
The perturbed electrostatic capacity   with respect to  infinity of  the conductor  $\Om_\ep$ is  given by
\begin{align}\label{capacity-Om-ep}
 \textrm{cap}(\Om_\ep)
=-\frac{1}{4\pi}\int_{\p\Om_\ep} \frac{\p u_\ep}{\p \tilde{n}}(y)d\tilde{\sigma}(y),
\end{align}
where  $\tilde{n}$  and $d\tilde{\sigma}$  are the unit outward normal  and the length element to the boundary $\p \Om_\ep$, respectively. Similarly to \eqref{asymptotic-infinity}, $u_\ep$ satisfies
\begin{align}\label{asymptotic-infinity-ep}
\ds u_\ep(x)=\frac{\textrm{cap}(\Om_\ep)}{ |x|} +O(\frac{1}{|x|^{2}})\quad \mbox{as } |x|\rightarrow \infty.
\end{align}

Our goal is to find  asymptotic expansions for the current, the electrostatic capacity, and the far-field of the electrostatic potential resulting from   small perturbations  on the surface of a conductor at equilibrium in free space. The main idea is to
 adopt the FE method to derive formal asymptotic expansions.
Then, based on layer potential techniques, we prove rigorously those asymptotic expansions.
 In connection with this, we refer to recent works in the  same  context \cite{KZ2, Z, Z2, LLZ, AKLZ1, AKLZ2, LTZ, KZ1, LZ}.

The  first achievement  of this paper is the following theorem,
a rigorous derivation of  the asymptotic expansion of the perturbed  current  $\p u_{\ep}/{\p \tilde {n}}$   on $\p \Om_\ep$ as $\ep \rightarrow 0$.
\begin{thm} \label{main-Theorem-current} Let $\tilde {x}=x+\ep h(x)n(x)\in \p \Om_{\ep}$, for $x\in \p \Om$. Let $u_\ep$ and $u$  be the solutions of \eqref{uep} and \eqref{u}, respectively.
 The following  asymptotic  expansion  holds:
\begin{align}
\ds  \pd{u_{\ep}}{\tilde{n}}(\tilde x)=\pd{u}{n}(x)+\ep \Big(2\tau(x) h(x) \pd{u}{n}(x)+\pd{v}{n}(x)\Big)+O(\ep^{2}),
\label{main-Theorem-equality-0}
\end{align}
where   the remainder $O(\ep^{2})$ depends  only
on the $\mathcal{C}^2$-norm of $X$ and  the $\mathcal{C}^1$-norm of $h$,    $\tau$  is the mean curvature of $\Om$, and  $v$  is  the unique solution to
\begin{equation}\label{v}
\left\{
  \begin{array}{ll}
\ds \Delta v=0 \quad\q  \mbox{ in } \RR^3 \backslash \overline{\Om}, \\
\nm \ds v=-h\pd{u}{n}  \quad \mbox{on  }  \p \Om,\\
\nm \ds  \lim_{|x|\rightarrow \infty}v(x) =0.
  \end{array}
\right.
\end{equation}
\end{thm}

The  second result of this paper is the following theorem, we rigorously derive   the asymptotic expansion of  $\textrm{cap}(\Om_\ep)$  as $\ep \rightarrow 0$.
\begin{thm} \label{main-Theorem}
 The following  asymptotic  expansion  holds:
\begin{align}
\ds  \textrm{cap}(\Omega_\ep)=\textrm{cap}(\Omega)+\frac{\ep}{4\pi}  \int_{\p \Omega}h\Big[\frac{\p u}{\p n }\Big]^2 d\sigma+O(\ep^2),
\label{main-Theorem-equality}
\end{align}
where   the remainder $O(\ep^2)$ depends  only
on the $\mathcal{C}^2$-norm of $X$ and the $\mathcal{C}^1$-norm of $h$.
\end{thm}
It is worth noticing that if $h$ has a constant sign  on $\p \Om$ that there exists $\ep_0 >0$ such that
for $\ep < \ep_0$, $\textrm{cap}(\Omega_\ep)-\textrm{cap}(\Omega)$ has the same sign as $h$.

The following theorem represents  the  third result of this paper,   a rigorous derivation of  the asymptotic expansion of the far-field of the  electrostatic potential $u_\ep$  as $\ep \rightarrow 0$.

\begin{thm}\label{main-theorem-3}
Let $u_\ep$ and $u$  be the solutions of \eqref{uep} and \eqref{u}, respectively. We have the following asymptotic expansion
 \begin{align} \label{far-field-expansion}
\ds u_\ep(x)=u(x)+\frac{\ep}{ 4 \pi |x|} \int_{\p \Omega}h\Big[\frac{\p u}{\p n }\Big]^2 d\sigma+O\big(\frac{\ep^2}{|x|}\big) +O\big(\frac{\ep}{|x|^{2}}\big)
\end{align}
as $|x|\rightarrow \infty$  and $ \ep \rightarrow 0$, where   the remainders $O({\ep^2}/{|x|})$  and  $O({\ep}/{|x|^{2}})$ depend  only
on the $\mathcal{C}^2$-norm of $X$,  the $\mathcal{C}^1$-norm of $h$, and  the $distance (origin, \Om)$.
\end{thm}

 These asymptotic expansions had not been established before this work.  They can be taken into account   in the  design of  conductors
to avoid negative effects due to small changes in their shapes. Our asymptotic expansions are still valid in the case of small perturbations
of one of the  walls  of a   condenser in electrostatic equilibrium, the FE method works well,
but more elaborate arguments are needed for the layer potential techniques  method.

The asymptotic  expansions can be used   to design effective algorithms to recover certain properties of the perturbation of the shape of an isolated conductor. That is, we would like to find a method for determining the shape of a  conductor  by taking  one or  a combination   of   current,  electrostatic capacity, and  electrostatic potential measurements.   One of solutions is to extend that optimization approach in  \cite{ABFKL} by using  electrostatic capacity  measurements.

This article is organized as follows. In the next section we  formally derive  the  asymptotic expansions  \eqref{main-Theorem-equality-0}, \eqref{main-Theorem-equality}, and  \eqref{far-field-expansion} by using the  FE method. In the section $2$, based on layer potential techniques  method we rigorously prove those in fact the formal expansions hold. In the last section, we derive an asymptotic  expansion of the   first eigenvector   of  the $L^2-$adjoint  of the   NP operator     resulting  from small perturbations  of the   surface   of its domain.

\section{Formal derivations: the FE method}
We refer to \cite{KZ2} for further   details on  the following concepts  and definitions.  Suppose $\p \Om$ has an orthogonal
parametrization  $X(\xi,\theta)$, that is,  there is an open subset  $\vartheta$  of $\RR^2$ such that
$ \p \Om:=\big\{x= X(\xi, \theta), (\xi, \theta)\in \vartheta\big\},$
 where  $X$ is a $\mathcal{C}^2$-function satisfying
 $(X_\xi:=\frac{d X}{d \xi})\cdot (X_\theta:=\frac{d
X}{d \theta})=0$ and $X_{\xi \theta}= X_{\theta \xi}$. We point out  that a revolution surface has an orthogonal parametrization.
 The vectors  $T_\xi=X_\xi/|X_\xi|$ and $T_\theta=X_\theta/|X_\theta|$ form an  orthonormal basis for the tangent plane to $\p \Om$ at $x=X(\xi, \theta)$. The tangential derivative on $\p \Om$ is defined  by $\pd{}{T}=\frac{\p }{\p T_\xi}T_{\xi}+\pd{}{T_\theta}T_{\theta}$.
 We denote by $\mathcal{{G}}$  the matrix of the first fundamental form with respect to the basis  $\{X_\xi, X_\theta\}$.

For $w\in \mathcal{C}^1(\vartheta)$, the gradient operator in local coordinates  satisfies
\begin{align}\label{local-gradient}
\ds \nabla_{\xi, \theta}w=\sqrt{\mathcal{G}_{11}}\pd{w}{T_\xi}T_\xi+\sqrt{\mathcal{G}_{22}}\pd{w}{T_\theta}T_\theta=\mathcal{G}^{\frac{1}{2}}\pd{w}{T}.
\end{align}
For $w\in \mathcal{C}^2(\vartheta)$, the  restriction of $\Delta $ in $\RR^3\backslash \p{\Om}$
to a neighbourhood of $\p \Om$ can be expressed as:
\begin{align}\label{local-Laplace}
\ds \Delta w=\frac{\p^2 w}{\p n^2}-2\tau \pd{w}{n}+\Delta_{\mathcal{G}}w,
\end{align}
where  $\Delta_{\mathcal{G}}$ is the Laplace–Beltrami operator  associated to $\mathcal{G}$ which is given by
\begin{align}\label{local-Betrami-Laplace}
\ds \Delta_{\mathcal{G}}w=\frac{1}{\sqrt{det\mathcal{G}}}
\nabla_{\xi, \theta}\cdot \Big(\sqrt{det\mathcal{G}} \mathcal{G}^{-1}\nabla_{\xi, \theta} w\Big).
\end{align}

We use $h(\xi, \theta)$ for simplifying the term $h(X(\xi, \theta))$
and $h_\xi(\xi, \theta)$, $h_\theta(\xi, \theta)$ for the tangential derivatives of $h(X(\xi, \theta))$.
Then,
$\tilde x =X(\xi, \theta)+\ep h(\xi, \theta) n(x)$ is a parametrization of $\p \Om_\ep$.

The following asymptotic expansions for the normal derivative $ \tilde{n}(\tilde {x})$
 and the length element $d\tilde{\sigma}(\tilde x)$ hold. For proofs,  see \cite{KZ2}.

\begin{lem} Let $\tilde {x}=x+\ep h(x)n(x)\in \p \Om_{\ep}$, for $x\in \p \Om$. Then,
the outward unit normal $\tilde{n}(\tilde{x})$ to $\partial
\Om_{\ep}$ at $\tilde{x}$ can be expanded uniformly as
\begin{equation}\label{nu-tilde}
\tilde{n}(\tilde{x}):=\frac{\tilde{X}_{\xi}\wedge
\tilde{X}_{\ta}}{|\tilde{X}_{\xi}\wedge \tilde{X}_{\ta}|}
=\sum_{k=0}^{\infty} \ep^k
 n_{k}(x),
 \end{equation}
where the vector-valued functions $n_{k}$ are uniformly bounded
regardless of $k$. In particular, for $x\in \partial \Om$, we have
\begin{align*}
 \ds n_{0}(x)=n(x),\q\q \q  n_{1}(x) = -\pd{h}{T}(x) T(x).
\end{align*}
 Likewise, the length element $d\tilde{\sigma}(\tilde x)$ has  the following
uniformly  expansion
\begin{align}\label{dsigma-ep}
\ds d\tilde{\sigma}(\tilde x):=|\tilde{X}_{\xi}\wedge \tilde{X}_{\ta}| d\xi
d\ta=\frac{|\tilde{X}_{\xi}\wedge
\tilde{X}_{\ta}|}{\sqrt{det(\mathcal{G})}}d\sigma(x)=\sum_{k=0}^{\infty}\ep
^k \sigma^{(k)}(x)d\sigma(x),
\end{align}
where $\sigma^{(k)}$ are uniformly bounded regardless of $k$ with
\begin{align*}
 \ds \sigma^{(0)}(x)=1, \q\q \q \sigma^{(1)}(x)=-2 h(x) \tau(x).
\end{align*}
\end{lem}

Let $u_\ep$ be the solution to \eqref{uep}. In order to derive a formal asymptotic expansion for $u_\ep$, we apply the FE method, see \cite{Z2,Coifman, KZ2, LLZ, Z}. Firstly, we expand $u_\ep$ in powers of $\ep$,
\begin{equation}\label{heatExp}
\ds u_\ep(x)=u_0(x)+\ep u_1(x)+\ep^2 u_2(x)+\cdots,\quad x\in
\RR^3\backslash \overline{\Om}_\ep,
\end{equation}
where  $u_l$ are  defined on $\RR^3\backslash \p {\Om}$. Since $u_l$ satisfy
\begin{equation}\label{u-l}
\left\{
  \begin{array}{ll}
\ds  \Delta u_{l} =0 \quad \mbox{in } \RR^3 \backslash \overline{\Om}, \\
\nm \ds \lim_{|x|\rightarrow \infty}u_{l}(x) =0.
  \end{array}
\right.
\end{equation}
In order to justify the first equation in \eqref{u-l}, we substitute \eqref{heatExp} in
the Laplace equation $\Delta u_\ep=0$ in $\RR^3\backslash \overline{\Om}_\ep$ to get  $\Delta u_{l} =0$
in $\RR^3\backslash \overline{\Om}_\ep$ for  $\ep>0$. Because $\ep$ is arbitrary,  we confirm that $\Delta u_{l}=0$
in $\RR^3\backslash \overline{\Om}$.

Let $\tilde x=x+\ep h(x) n(x) \in \p \Om_\ep$, for $x\in \p \Om$. The following Taylor expansion  holds:
\begin{align}\label{Taylor-normaluep-3D}
\ds \pd{u_\ep}{\tilde{n}}(\tilde x)=& \Big(\nabla u_0\big(x+\ep h(x) n(x)\big)+\ep \nabla u_1 \big(x+\ep h(x) n(x)\big)\Big)\cdot \tilde{n}(\tilde x)+O(\ep^2)\nonumber\\
\nm \ds =&\Big(\nabla u_0(x)+\ep h(x)\nabla^2 u_0(x)n(x)+\ep \nabla u_1( x)\Big)\cdot \Big(n(x)-\ep \pd{h}{T}(x)T(x)\Big)+O(\ep^2)\nonumber\\
\nm \ds  =& \pd{u_0}{n}(x)+\ep  h(x) \pdl{u_0}{n^2}(x)
+ \ep \pd{u_1}{n}(x)-\ep \pd{h}{T}(x)\cdot \pd{u_0}{T}(x)+O(\ep^2)\nonumber\\
\nm \ds  =& \pd{u_0}{n}(x)+2\ep \tau(x) h(x) \pd{u_0}{n}(x)
+ \ep \pd{u_1}{n}(x)\nonumber\\
\nm \ds  &-\ep \pd{h}{T}(x)\cdot \pd{u_0}{T}(x)-\ep h (x)\Delta_{\mathcal{G}} u_0(x)+O(\ep^2).
\end{align}
To justify the last equality, we use the representation of the Laplace operator on $\p \Om$ given in \eqref{local-Laplace}
\begin{align*}
\ds 0=\Delta {u_0}=\frac{\p^2 {u_0}}{\p n^2}-2\tau \pd{{u_0}}{n}+\Delta_{\mathcal{G}} u_0 \q \mbox{on } \p \Om.
\end{align*}
For $\tilde x=x+\ep h(x) n(x) \in \p \Om_\ep$. We have the following  Taylor expansion
\begin{align}\label{Taylor-uep-3D}
\ds  u_\ep(\tilde x) &=u_0\big(x+\ep h(x) n(x)\big)+ \ep u_1\big(x+\ep h(x) n(x)\big)+O(\ep^2)\nonumber\\
\nm \ds &= u_0(x)+\ep h(x) \pd{u_0}{n}(x)+ \ep u_1(x)+O(\ep^2).
\end{align}
Using   the boundary condition $u_\ep(\tilde x)=1$ for $\tilde x =\Psi(x) \in \p \Om_\ep$,  we   obtain  from \eqref{Taylor-uep-3D} that
\begin{align}\label{boundary-conditions-u0-u1}
\ds u_0(x)=1,\q  u_1(x)=- h(x) \frac{\p u_0}{\p n}(x), \q\mbox{for }x \in \p \Om.
\end{align}
By the uniqueness of the  PDE problems   \eqref{u} and \eqref{v},   we get  $u_0\equiv u$  and $u_1 \equiv v$ in $\RR^3 \backslash \overline{\Om}$. The  fourth and fifth terms  in \eqref{Taylor-normaluep-3D} vanish, this is because $u=1$ on $\p \Om$ which implies that ${ \p u}/{\p T}=0$ on $\p \Om$.
Then Theorem  \ref{main-Theorem-current} immediately follows from \eqref{Taylor-normaluep-3D} formally.

A change of variables $y=\Psi_{\ep}(x)$   for $x\in \p \Om$ in \eqref{capacity-Om-ep} gives
\begin{align}\label{7}
\ds \textrm{cap}(\Om_\ep)=-\frac{1}{4\pi}\int_{\p\Om} \frac{\p u_\ep}{\p \tilde{n} }(\tilde x) u_\ep(\tilde x) d\tilde{\sigma}(\tilde x).
\end{align}
It  then follows from   \eqref{dsigma-ep}, \eqnref{Taylor-normaluep-3D}, and \eqref{Taylor-uep-3D} that
\begin{align}\label{1234}
\ds  \textrm{cap}(\Om_\ep)=&-\frac{1}{4\pi}\int_{\p\Om} \frac{\p u}{\p n } u d\sigma
-\frac{\ep}{4\pi} \int_{\p \Om}h \Big [\frac{\p u}{\p n } \Big]^2d\sigma\nonumber\\
\nm\ds &-\frac{\ep}{4\pi} \int_{\p \Om} \frac{\p v}{\p n } u d\sigma
-\frac{\ep}{4\pi} \int_{\p \Om} \frac{\p u}{\p n } v d\sigma+O(\ep^2).
\end{align}
By Green's identity and \eqref{u-l}, we have
\begin{align*}
\ds  \int_{\p \Om} \frac{\p v}{\p n } u d\sigma
=\int_{\p \Om} \frac{\p u}{\p n } v d\sigma.
\end{align*}
We get from  \eqref{boundary-conditions-u0-u1}   that
\begin{align}\label{333}
\int_{\p \Om} \frac{\p v}{\p n } u d\sigma
+\int_{\p \Om} \frac{\p u}{\p n } v d\sigma=2\int_{\p \Om} \frac{\p u}{\p n } v d\sigma=-2\int_{\p \Om}h \Big [\frac{\p u}{\p n } \Big]^2d\sigma.
\end{align}
Thus, by     \eqref{capacity-Om}, \eqref{1234},  and \eqref{333}, we formally obtain the desired Theorem \ref{main-Theorem}, $i.e.$,
\begin{align}\label{formal-cap}
\ds \textrm{cap}(\Omega_\ep)=\textrm{cap}(\Omega)+\frac{\ep}{4\pi} \int_{\p \Omega}h\Big[\frac{\p u}{\p n }\Big]^2 d\sigma+O(\ep^2).
\end{align}
As a direct consequence of  \eqref{asymptotic-infinity}, \eqref{asymptotic-infinity-ep}, and \eqref{formal-cap}, the leading order term in the asymptotic expansion   of the far-field $u_\ep-u$  in Theorem \ref{main-theorem-3} holds formally
 \begin{align}\label{far-field-formal}
\ds u_\ep(x)-u(x)=\frac{\ep}{ 4 \pi |x|} \int_{\p \Omega}h\Big[\frac{\p u}{\p n }\Big]^2 d\sigma+O(\frac{\ep^2}{|x|}) +O(\frac{1}{|x|^{2}})
\end{align}
  as $\ep \rightarrow 0$ and  $\ep>>1/|x|$. By the layer potential techniques  method  we will  prove in the subsection \ref{proof-theorem-3} the  asymptotic expansion \eqref{far-field-formal} with  a remainder   $O({\ep^2}/{|x|}) +O({\ep }/{|x|^{2}})$  as $\ep\rightarrow 0$ and  $|x|\rightarrow \infty$ which is more better than $O({\ep^2}/{|x|}) +O({1}/{|x|^{2}})$ as $\ep \rightarrow 0$ and  $\ep>>1/|x|$.
\section{Layer potential techniques  method}
\subsection{ Definitions and Preliminary results}
Let $\Om$ be a bounded $\mathcal{C}^2$-domain. Let $\Gamma(x)$  be
the fundamental solution of the Laplacian $\Delta$ in $\RR^3$:
$\Gamma(x)=-\frac{1}{4\pi|x|}$. The single  and double layer
potentials of the density function $\phi$ on $\p \Om$ are  defined by
 \begin{align}
\ds  \Scal_{\Om} [\phi](x)
 & = \int_{\partial \Om} \Gamma(x-y) \phi(y)d\sigma(y), \quad x \in
 \RR^3,\label{single-layer}\\
\nm \ds  \Dcal_{\Om} [\phi] (x)
 & = \int_{\partial \Om} \frac{\p }{\p n(y)}\Gamma(x-y) \phi(y)d\sigma(y), \quad x \in \RR^3 \setminus \p
 \Om.\label{double-layer}
\end{align}
We note  that for $x\in \RR^3\backslash \p \Om$ and $y\in \p \Om$, $\Gamma(x-y)$ and $\frac{\p}{\p n(y)}\Gamma(x-y)$
are $L^{\infty}$-functions in $y$  and harmonic in $x$  and their  behaviors  when $|x|\rightarrow +\infty  $  are given by
\begin{align}\label{asymptotic-green-function-infinity}
\ds \Gamma(x-y)=O(\frac{1}{|x|}), \q\q \q\q \frac{\p}{{\p n(y)}} \Gamma(x-y)=O(\frac{1}{|x|^2}).
\end{align}
Therefore, we readily see that $\Dcal_{\Om}[\phi]$ and $\Scal_{\Om}[\phi]$ are well defined and harmonic in $\RR^3\backslash \p \Om$ and  satisfy
\begin{align}\label{single-double-infity}
\ds \Scal_{\Om} [\phi](x)=O(\frac{1}{|x|}), \q \q\q \Dcal_{\Om} [\phi] (x)=O(\frac{1}{|x|^2}), \q \mbox{as } |x|\rightarrow +\infty.
\end{align}

 We denote  $\pd{w}{n}|_{\pm}=n\cdot \nabla w^{\pm}|_{\p \Om}$, where $w^{+}=w|_{\RR^3\backslash \overline{\Om}}$ and  $w^{-}=w|_{ \Om}$.   The following formulae give the jump relations obeyed by the double
layer potential and by the normal derivative of the single layer
potential. For proofs, see  \cite{Folland76, Book-Ammari}.
\begin{align}
\ds \Scal_{\Om} [\phi] \big |_{+}(x) & = \Scal_{\Om} [\phi] \big |_{-}(x) \quad
\mbox{a.e. } x \in \p {\Om} \label{Single-boundary}, \\
\nm \ds \pd{(\Scal_{\Om} [\phi])}{T} \Big |_{+}(x) & =  \pd{(\Scal_{\Om} [\phi])}{T} \Big |_{-}(x)\quad \mbox{a.e. } x \in \p {\Om} \label{nuS}, \\
\nm \ds \pd{(\Scal_{\Om} [\phi])}{n} \Big |_{\pm}(x) & = \big (\pm \frac{1}{2} I
+ (\Kcal_{\Om})^* \big ) [\phi] (x) \quad
\mbox{a.e. } x \in \p {\Om} \label{nuS}, \\
\nm \ds  (\Dcal_{\Om} [\phi]) \big|_{\pm} (x) & = \big(\mp \frac{1}{2}
I + \Kcal_{\Om} \big) [\phi](x) \quad \mbox{a.e. } x \in \p {\Om}, \label{doublejump-h}
\end{align}
for $\phi \in L^2(\p {\Om})$, where $\Kcal_{\Om}$ is the NP  operator defined by
$$ \Kcal_{\Om} [\phi] (x) =\frac{1}{4\pi} \int_{\partial {\Om}} \frac{\la
 y-x,n(y)\ra}{|x-y|^3}\phi(y)d\sigma(y), $$ and
$\Kcal_{\Om}^*$ is the $L^2$-adjoint operator of the NP  operator $\Kcal_{\Om}$, that is,
\begin{equation}\label{KD*}
 \Kcal_{\Om}^*
[\phi] (x) =\frac{1}{4\pi} \int_{\partial {\Om}} \frac{\la
 x-y,n(x)\ra}{|x-y|^3}\phi(y)d\sigma(y).
 \end{equation}
 The operators $\Kcal_\Om$ and
$\Kcal_\Om^*$ are singular integral operators and bounded on $L^2(\p
\Om)$. Because $\Om$ has a $\mathcal{C}^2$ boundary,  $\pd{(\Dcal_\Om
[\phi])}{n}$ does not have a jump across $\p \Om$, that is,
\begin{align}\label{normal-D}
\ds \pd{(\Dcal_{\Om} [\phi])}{n} \Big |_{+}(x)=\pd{(\Dcal_{\Om} [\phi])}{n} \Big
|_{-}(x),\quad x\in \p {\Om}.
\end{align}

Let
$
\ds  W^2_1 (\p {\Om}) : = \{ f \in L^2(\p {\Om}) : { \p f}/{\p T} \in L^2(\p {\Om})
\}.
$
 The following lemma is of importance to us. For proof, see for example
\cite{book-1}.

\begin{lem}\label{invertibity-single-layer} Let $\Om$  be a bounded Lipschitz  domain  in $\RR^3$. Then
\begin{enumerate}[(i)]
\item $\Scal_{\Om}: L^2 (\p {\Om}) \longrightarrow W^2_1 (\p {\Om})$  has a bounded inverse.
\item $\Kcal_{\Om}: W^2_1 (\p {\Om}) \longrightarrow W^2_1 (\p {\Om})$ is a bounded operator.
\end{enumerate}
\end{lem}

We will need the following lemma which was
obtained in  \cite{Folland76}; see also  \cite{Book-Ammari}.

\begin{lem}\label{double-layer potential- behavior boundary}
If $\Om$ is a bounded $\mathcal{C}^2$-domain. Then $\Dcal_{\Om}[1](x)=0  $  for $x\in \RR^3 \setminus
 \overline{\Om}$,  $\Dcal_{\Om}[1](x)=1  $  for $x\in
{\Om}$,  and  $ \Kcal_{\Om}[1]=\frac{1}{2}$  for $x\in \p \Om$.
\end{lem}

\subsection{Asymptotic of layer potentials}
By using the change of variable  $z=\Psi_{\ep}(y)=\tilde{y}$   for  $y\in \p \Om$  and  $z\in \p \Om_\ep$, we write
\begin{align*}
\ds  \Scal_{{\Om_\ep}}[\tilde{\psi}](\tilde{x})=-\frac{1}{   4\pi}
  \int_{\partial {\Om}_\epsilon}\frac{1}
 {|\tilde{x}-z|}\tilde
 {\psi}(z)d\tilde{\sigma}(z) =-\frac{1}{   4\pi}
  \int_{\partial {\Om}}\frac{1}
 {|\tilde{x}-\tilde{y}|}\tilde
 {\psi}(\tilde{y})d\tilde{\sigma}(\tilde{y}),\quad  \tilde{x} \in \p \Om_\ep,
 \end{align*}
for any density $\tilde \psi\in L^2(\partial \Om_\ep)$.
For $(\xi,\ta),(\alpha,\beta)\in \vartheta.$ Set
$$
x=X(\xi,\ta),\quad\quad \tilde x=\tilde X(\xi,\ta)=x+\ep
h(\xi,\ta)n(x),
$$
$$
y=X(\alpha,\beta),\quad\quad \tilde y=\tilde X(\alpha,\beta)=y+\ep
h(\alpha,\beta)n(y),
$$
and hence
\begin{equation*}
\tilde x-\tilde y= x- y+\ep
\Big(h(\xi,\ta)n(x)-h(\alpha,\beta)n(y)\Big).
\end{equation*}
This gives
\begin{align}\label{01}
\ds |\tilde x-\tilde y|
=&|
x-y|\left(1+2\ep \frac{\la
x-y,h(x)n(x)-h(y)n(y)\ra}{|
x-y|^2}+\ep^2\frac{\big|h(x)n(x)-h(y)n(y)\big|^2}{|x-y|^2}\right)^{\frac{1}{2}}\nonumber\\
:=&| x-y|\Big(1+2\ep F(x,y)+\ep^2G(x,y)\Big)^{\frac{1}{2}}.
\end{align}
We have $h\nu \in \mathcal{C}^1(\p \Om)$. Then,  one can easily see that
 $$
 |F(x,y)|+|G(x,y)|^{\frac{1}{2}} \leq C \| X \|_{\mathcal{C}^2{(\p \Om)}} \|h \|_{\mathcal{C}^1(\p \Om)}\q \mbox{for } x,y\in \p \Om.
 $$
Therefore, it follows from \eqref{dsigma-ep} and \eqref{01} that
\begin{align}\label{norm}
\ds \frac{1}{|\tilde x-\tilde y|}d \tilde{\sigma}(\tilde y)=&\frac{1}{| x-y|}\Big(1+2\ep
F(x,y)+\ep^2G(x,y)\Big)^{-\frac{1}{2}} \times \Big(\sum_{k=0}^{\infty}  \ep^k \sigma_k(y) d\sigma(y)\Big)\nonumber \\
 :=&\frac{1}{| x-y|} \sum_{k=0}^{\infty}  \ep^k \mathbb{L}_{k}(x,y) d\sigma(y),
\end{align}
where
$
 |\mathbb{L}_{k}(x,y)| \leq C \| X \|_{\mathcal{C}^2{(\p \Om)}} \|h \|_{\mathcal{C}^1(\p \Om)}, \mbox{ for } x,y\in \p \Om.
$
In particular
\begin{align*}
\ds \mathbb{L}_{0}(x,y)=1,\q\q\q  \mathbb{L}_{1}(x,y)=- F(x,y)-2 \tau(y) h(y).
\end{align*}
Introduce a sequence of integral operators
$(\Scal_{\Om}^{(k)})_{k\in\NN}$, defined for any  $\psi \in L^2(\p \Omega)$
by
\begin{align*}
 \ds \Scal_{\Om}^{(k)}\psi(x):=-\frac{1}{4\pi}\int_{\p
\Om } \frac{\mathbb{L}_k(x,y)}{| x-y|}\psi(y)d\sigma(y) \quad \mbox{ for } k\geq 0.
\end{align*}
Note that $\Scal_{\Om}^{(0)}=\Scal_{\Om}$ and
\begin{align*}
\ds \Scal_{\Om}^{(1)}[\psi](x)=&\frac{1}{2 \pi}\int_{\p
\Om }\frac{1}{| x-y|}\tau(y)h(y)\psi(y) d\sigma(y)\nonumber\\
\nm \ds &+  \frac{h(x)}{4\pi}\int_{\p
\Om } \frac{\la x-y,n(x)\ra
}{|x-y|^3}\psi(y) d\sigma(y)-\frac{1}{4 \pi} \int_{\p
\Om } \frac{\la x-y,n(y)\ra
}{|x-y|^3}h(y)\psi(y) d\sigma(y)\nonumber\\
\nm \ds =&- 2\Scal_{\Om}[\tau h \psi](x)+h(x) \frac{\p(\Scal_{\Om}[\psi])}{\p n}\Big |_{\pm}(x)+ \Dcal_{\Om}[h\psi]\big |_{\pm}(x) \quad  \mbox{for  } x\in \p \Om.
\end{align*}
It is easily to prove that the operator $\Scal_{\Om}^{(k)}$ with the
kernel $-\frac{1}{4 \pi} {\mathbb{L}_k(x,y)}/{| x-y|}$  is bounded on $L^2(\p \Om)$. See \cite[Proposition 3.10]{Folland76}.

Let  $\tilde{x}=\Psi_\ep(x)=x+\ep  h(\xi,\ta) n(x)$ for
$x=X(\xi,\ta)\in \p\Om$. The following estimate  holds:
\begin{align*}
\ds \Big\|\Scal_{{\Om}_\ep}[\tilde{{\psi}}]\circ \Psi_\ep-\Scal_{\Om}[\psi]-\sum_{k=1}^{N}\ep^k \Scal_{\Om}^{(k)}[\psi]\Big\|_{L^{2}(\p {\Om})}\leq C \ep^{N+1}\big\|\psi \big\|_{L^2(\p {\Om})},
\end{align*}
where $\psi:=\tilde{\psi}\circ \Psi_\ep$ and $C$ depends
only on $N$, $\|X\|_{\mathcal{C}^2(\p\Om)}$, and
$\|h\|_{\mathcal{C}^1(\p\Om)}$.

We have
\begin{align*}
\ds \nabla_x \frac{1}{|\tilde x-\tilde y|}\cdot T(x) d\tilde{\sigma}(\tilde y)=& \sum_{k=0}^{\infty}  \ep^k \Big[\frac{\la x-y, T(x)\ra}{| x-y|^3}\mathbb{L}_{k}(x,y) + \frac{\la \nabla_x \mathbb{L}_{k}(x,y), T(x)\ra}{| x-y|}\Big]d\sigma(y)\\
\nm \ds :=& \sum_{k=0}^{\infty}\ep^k \mathbb{K}_{k}(x,y)d\sigma(y).
\end{align*}
By looking at the $\nabla F(x,y) \cdot T(x)$ and $\nabla G(x,y) \cdot T(x)$,  we confirm that  $\mathbb{K}_{k}(x,y)$
is a combination linear as following
\begin{align*}
\ds \mathbb{K}_{k}(x,y)=&\alpha_k(x,y)\frac{\la x-y, T(x)\ra}{| x-y|^3}+\beta_k(x,y)\frac{\la h(x)n(x)-h(y)n(y), T(x)\ra}{| x-y|^3}\\
\nm \ds &+\gamma_k(x,y)\frac{\la x-y, n(x)\ra}{| x-y|^3}+\lambda_k(x,y)\frac{\la h(x)n(x)-h(y)n(y), n(x)\ra}{| x-y|^3},
\end{align*}
where $|\alpha_k(x,y)|+|\beta_k(x,y)|+|\gamma_k(x,y)|+|\lambda_k(x,y)|\leq C \| X \|_{\mathcal{C}^2{(\p \Om)}} \|h \|_{\mathcal{C}^1(\p \Om)}, \mbox{ for } x,y\in \p \Om.$

It is easily to prove that the operator $\p\Scal_{\Om}^{(k)}/{\p T}$ with the
kernel $-\frac{1}{4 \pi} \mathbb{K}_{k}(x,y)$  is bounded on $L^2(\p \Om)$. In fact, it is an
immediate consequence of the celebrate theorem of
Coifman-McIntosh-Meyer,  see \cite{CMM82}. Therefore, the following estimate holds:
\begin{align*}
\ds \bigg\|\pd{\Scal_{{\Om}_\ep}[{\tilde{\psi}}]\circ \Psi_\ep}{T}-\pd{\Scal_{\Om}[\psi]}{T}-\sum_{k=1}^{N}\ep^k \pd{\Scal_{\Om}^{(k)}[\psi]}{T}\bigg\|_{L^{2}(\p {\Om})}\leq C \ep^{N+1}\big\|\psi \big\|_{L^2(\p {\Om})}.
\end{align*}

 The result of the above asymptotic analysis is summarized in the following theorem.
\begin{thm}\label{asymptotic-single-layer}
There exists  $C$ depending only on  $\|X\|_{\mathcal{C}^2(\p \Om)}$ and $\|h\|_{\mathcal{C}^1(\p \Om)}$, such that for any $\tilde{{\psi}} \in L^2(\p \Om_\ep)$,  we have
\begin{align}\label{asymptotic-single-layer-bounded}
\ds \Big\|\Scal_{{\Om}_\ep}[\tilde{\psi}]\circ \Psi_\ep-\Scal_{\Om}[\psi]-\sum_{k=1}^{N}\ep^k \Scal_{\Om}^{(k)}[\psi]\Big\|_{W^{2}_{1}(\p {\Om})}\leq C \ep^{N+1}\big\|\psi \big\|_{L^2(\p {\Om})},
\end{align}
where $\psi:=\tilde{\psi}\circ \Psi_\ep$.
\end{thm}

For $\psi \in L^2(\p \Om)$, we introduce
\begin{align*}
\ds \Kcal_{\Om}^{(1)}[\psi](x)
=&2 \left (\tau h
\frac{\p(\Scal_{\Om}[\psi])}{\p n }-\frac{\p(\Scal_{\Om}[\tau h \psi])}{\p n}\right)\Big| _{\pm}(x)+\frac{\p(\Dcal_{\Om}[h\psi])}{\p n}(x)\nonumber\\
\nm \ds &-\frac{1}{\sqrt{det(\mathcal{G}})}\left(\nabla_{\xi,\ta} \cdot
\Big(h\sqrt{det(\mathcal{G})}\mathcal{G}^{-1}\nabla_{\xi,\ta}\Scal_{\Om}[\psi]\Big)\right)(x), \quad  \mbox{for  } x\in \p \Om.
\end{align*}
It was proved in \cite{KZ2} that the operator   $ \Kcal_{\Om}^{(1)}$  is bounded in $L^2(\p {\Om})$ and  the following proposition  holds.
\begin{prop}\label{asymptotic-normal-single-layer}
There exists  $C$ depending only on  $\|X\|_{\mathcal{C}^2(\p \Om)}$ and $\|h\|_{\mathcal{C}^1(\p \Om)}$, such that for any $\tilde{\psi} \in L^2(\p {\Om}_\ep)$,  we have
\begin{align}
\bigg\|\pd{\Scal_{{\Om}_\ep}[\tilde{\psi}]}{\tilde{n}}\circ \Psi_\ep\Big|_{\pm}-\pd{\Scal_{\Om}[\psi]}{n}\Big|_{\pm}-\ep \Kcal_{\Om}^{(1)}[\psi]\bigg\|_{L^2(\p {\Om})}\leq C \ep^2\big\|\psi \big\|_{L^2(\p {\Om})},
\end{align}
where $\psi:=\tilde{\psi}\circ \Psi_\ep$.
\end{prop}

\subsection{Proofs of Theorems }
The following lemma  is of use to us.
\begin{lem}\label{important-lemma} Let $f \in W^{2}_{1}(\p \Om)$. The solution  of the following problem
\begin{equation}\label{U}
\left\{
  \begin{array}{ll}
\ds \Delta w=0 \quad \mbox{ in } \RR^3 \backslash \overline{\Om}, \\
\nm \ds w=f \q \quad \mbox{on  }  \p \Om,\\
\nm \ds  \lim_{|x|\rightarrow \infty}w(x) =0,
  \end{array}
\right.
\end{equation}
is represented as
\begin{align}\label{representation-formula-U}
w(x)=\Scal_{\Om}[\phi](x)-\Dcal_{\Om}[f](x),\q x\in \RR^3 \backslash \overline{\Om},\q \phi:=\pd{w}{n}\Big|_{\p \Om},
\end{align}
where $\phi \in L^2(\p \Om)$ satisfies  the following integral equation
\begin{align*}
\Scal_{\Om}[\phi]=(\frac{1}{2}I+\Kcal_{\Om})[f]\q \mbox{on } \p\Om.
\end{align*}
The representation formula \eqref{representation-formula-U} is unique.
\end{lem}
\proof
Consider the following problem
\begin{equation}\label{w}
\left\{
  \begin{array}{ll}
\ds \Delta U=0 &\mbox{ in } \RR^3 \backslash \p \Om, \\
\nm \ds U|_{+}-U|_{-}=f & \mbox{on  }  \p \Om,\\
\nm \ds \pd{U}{n}\Big |_{+}-\pd{U}{n}\Big |_{-}=\phi & \mbox{on  }  \p \Om,\\
\nm \ds  U(x) =O(1/|x|)& \mbox{as  } x\rightarrow \infty,
  \end{array}
\right.
\end{equation}
Let $U_1=\Scal_{\Om}[\phi]-\Dcal_{\Om}[f]$ in $\RR^3$. It follows from \eqref{single-double-infity} that $U_1(x)=O(1/|x|)$
and hence $U_1$ is a solution  of \eqref{w}  by the jump formulae \eqref{Single-boundary}--\eqref{doublejump-h}, and \eqref{normal-D}.
If we put $U_2=w$ in $\RR^3\backslash \overline{\Om}$ and $U_2\equiv 0$  in $\Om$, then $U_2$ is also a solution of \eqref{w}.
Therefore, in order to prove \eqref{representation-formula-U}, it suffices to show that the problem \eqref{w} has a unique solution  in
$W_{\mbox{loc}}^{1,2}(\RR^3\backslash \p \Om)$.

Suppose that $U\in W_{\mbox{loc}}^{1,2}(\RR^3\backslash \p \Om)$ is a solution  of \eqref{w}  with $f=\phi=0$. Then  $U$ is the weak solution of $\Delta U=0$ in the entire domain $\RR^3$. Therefore, for a  large $R$,
$$
\int_{B_R(0)} |\nabla U|^2 dx = \int_{\p B_R(0)} U \pd{U}{n} d\sigma(x)= - \int_{\RR^3\backslash \overline{B_R(0)}} |\nabla U|^2 dx\leq 0,
$$
where $B_{R}(0)=\{|x|< R\}$. This inequality holds for all $R$ and hence $U$ is constant. Since $U(x) \rightarrow 0$ at the infinity,
 we conclude that $U=0$.

From\eqref{representation-formula-U} and \eqref{doublejump-h}, we get $\Scal_{\Om}[\phi]=(\frac{1}{2}I+\Kcal_{\Om})[f] \mbox{ on } \p\Om.$ For $f \in W^{2}_{1}(\p \Om)$, we have $\Kcal_{\Om}[f]\in W^{2}_{1}(\p \Om)$ and hence    $(\frac{1}{2}I+\Kcal_{\Om})[f]\in W^{2}_{1}(\p \Om)$. It then follows from Lemma \ref{invertibity-single-layer} $(i)$ that $\phi $ is unique and belongs to  $L^2(\p \Om)$. Therefore,  the  representation formula \eqref{representation-formula-U} is unique.

According to Lemmas \ref{double-layer potential- behavior boundary} and \ref{important-lemma},
 the solution  $u$    to \eqref{u} has  the following representation formula
\begin{align}\label{representation-u}
\ds u(x)=\Scal_{\Om}[\phi_0](x)-\Dcal_{\Om}[1](x)=\Scal_{\Om}[\phi_0](x), \q x \mbox{ in }\RR^3\backslash \overline{\Om},\q \phi_0:= \pd{u}{n}\Big|_{\p \Om},
\end{align}
where $\phi_0$  satisfies the integral equation
\begin{align}\label{phi}
\Scal_{\Om}[\phi_0]=1 \q \mbox{on }\p \Om.
\end{align}
 Similarly to \eqref{representation-u}, the solution $u_\ep$
to \eqref{uep}  is represented by:
\begin{align}\label{representation-u-ep}
\ds u_\ep(x)=\Scal_{{\Om}_\ep}[\phi_\ep](x), \q x \mbox{ in }\RR^3\backslash \overline{{\Om}}_\ep,
\end{align}
where  $\phi_\ep$ is the unique solution to
\begin{align}\label{phi-ep}
\Scal_{\Om_\ep}[\phi_\ep]=1  \q \mbox{on }\p \Om_\ep.
\end{align}
The following Lemma holds.

\begin{lem}\label{estimes-Energy-Lemma}
There exists $C$ depending only on the $\mathcal{C}^2$-norm of $X$ and $\mathcal{C}^1$-norm of $h$ such that
\begin{align}\label{estimation-energy-1}
\ds \big \|\phi_\ep \circ \Psi_\ep -\phi_0  \big \|_{L^{2}(\p \Om)}\leq C \ep,
\end{align}
where $\phi_0 $ and $\phi_\ep$ are defined in \eqref{phi} and  \eqref{phi-ep}, respectively.
\end{lem}
\proof
Let $x\in \p\Om$, then $\tilde x =\Psi_\ep(x)=x+\ep h(x) n(x) \in \p \Om_\ep$. According to \eqref{phi} and  \eqref{phi-ep} we have
\begin{align}\label{Boundary-equality-zero}
\ds \Scal_{\Om_\ep}[\phi_\ep]\circ \Psi_\ep(x)=\Scal_{\Om}[\phi_0](x),\q x\in \p \Om.
\end{align}
It then follows from Theorem \ref{asymptotic-single-layer} that
$$
\Scal_{\Om}[\phi_\ep\circ \Psi_{\ep}-\phi_0](x)=O(\ep),\q x\in \p \Om,
$$
with $O(\ep)$ is bounded in $W_{1}^2(\p  \Om)$ by $C\ep$ for some constant $C>0$ depending only on the $\mathcal{C}^2$-norm of $X$ and $\mathcal{C}^1$-norm of $h$.  Clearly the  desired estimate \eqref{estimation-energy-1} immediately follows from Lemma \ref{invertibity-single-layer}  $(i)$.

\subsubsection{Proof of Theorem \ref{main-Theorem-current}}\label{proof-theorem-1}
Let $\phi_\ep$ and $\phi_0$ be the solutions of the integral equations \eqref{phi-ep} and \eqref{phi}, respectively. We denote by $\phi:=\phi_\ep\circ \Psi_\ep$. Thanks to Lemma \ref{estimes-Energy-Lemma}, we write
\begin{align}\label{decomposition-phi-ep}
\phi= \phi_0 +\ep \phi_1,
\end{align}
with $\phi_1$ is bounded in $L^2(\p \Om)$ and still  depends on $\ep$. Define
\begin{align}\label{}
\ds v_\ep(x)=\Scal_{ \Om} [\phi_1](x)-2 \Scal_{ \Om}[\tau h \phi_0](x)+\Dcal_{\Om}[h\phi_0](x), \quad x \in \RR^3\backslash \overline{\Om}.
\end{align}
It follows from Proposition \ref{asymptotic-normal-single-layer} and \eqref{decomposition-phi-ep} that
\begin{align}\label{8}
\ds \Big\| \pd{u_{\ep}}{\tilde{n}}\circ \Psi_\ep-\pd{u}{n}-2\ep\tau h \pd{u}{n}-\ep\pd{v_\ep}{n}\Big\|_{L^2(\p {\Om})}\leq C \ep^2\big\|\phi \big\|_{L^2(\p {\Om})}.
\end{align}
Turning to Theorem \ref{asymptotic-single-layer} and \eqref{Boundary-equality-zero}, we confirm that
\begin{align}\label{4444}
\ds \Big\|\Scal_{{\Om}}[{\phi}_0]-\Scal_{\Om}[\phi]-\ep \Scal_{\Om}^{(1)}[\phi]\Big\|_{W^2_{1}(\p \Om)}\leq C \ep^2\big\|\phi \big\|_{L^2(\p {\Om})}.
\end{align}
Substitute $\phi= \phi_0 +\ep \phi_1$  in \eqref{4444}, we get
\begin{align}\label{444}
\ds \Big\|h \pd{\Scal_{{\Om}}[{\phi}_0]}{n}\Big|_{+}+\Scal_{\Om}[\phi_1]-2\Scal_{\Om}[\tau h\phi_0]+\Dcal_{\Om}[ h\phi_0]\big|_{+}
\Big\|_{W^2_{1}(\p \Om)}\leq C \ep\big\|\phi \big\|_{L^2(\p {\Om})},
\end{align}
that is,
\begin{align}\label{99}
\ds \|v_\ep -v \|_{W^2_{1}(\p \Om)} \leq C \ep.
\end{align}
Since  $v_\ep-v$  is harmonic  in $\RR^3 \backslash \overline{\Om}$, we obtain  from Lemma \ref{important-lemma} that
\begin{align*}
\ds (v_\ep- v)(x)=\Scal_{ \Om} \Big[\pd{v_\ep}{n}-\pd{v}{n}\Big](x)-\Dcal_{ \Om} [v_\ep- v](x), \q  x\in  \RR^3\backslash \overline{\Om},
\end{align*}
and therefore, we deduce from \eqref{doublejump-h} that
\begin{align}\label{9}
\ds \Scal_{ \Om} \Big[\pd{v_\ep}{n}-\pd{v}{n}\Big]=
\big(\frac{1}{2}+\Kcal_{\Om}\big)[v_\ep-v] \q \mbox{on } \p \Om.
\end{align}
It then follows from Lemma \ref{invertibity-single-layer},  \eqref{99},  and \eqref{9}  that
\begin{align}\label{Estimation-normal-v-ep-minus-normal-v}
\ds \Big \|\pd{v_\ep}{n}-\pd{v}{n}  \Big\|_{L^2(\p \Om)} \leq C \ep.
\end{align}
Finally, we prove the theorem \ref{main-Theorem-current}  as desired from \eqref{8} and \eqref{Estimation-normal-v-ep-minus-normal-v}.

\subsubsection{Proof of Theorem \ref{main-Theorem}}\label{proof-theorem-2}
It follows from  Theorem \ref{asymptotic-single-layer},
Proposition \ref{asymptotic-normal-single-layer}, \eqref{dsigma-ep}, and  \eqref{7} that
\begin{align*}
\ds \textrm{cap}(\Om_\ep)=&-\frac{1}{4\pi} \int_{\p \Om}\pd{\Scal_{\Om_\ep}[\phi_\ep]}{\tilde{n}}\Big|_{+}(\tilde x) \Scal_{\Om_\ep}[\phi_\ep](\tilde x)d\tilde{\sigma}(\tilde x)\\
\nm\ds =&-\frac{1}{4\pi} \int_{\p \Om}\pd{\Scal_{\Om}[\phi]}{n}\Big|_{+} \Scal_{\Om}[\phi]d\sigma
-\frac{\ep}{4\pi} \int_{\p \Om}\bigg[-2\pd{\Scal_{\Om}[\tau h\phi]}{n}\Big|_{+}+\pd{\Dcal_{\Om}[ h\phi]}{n} \bigg]\Scal_{\Om}[\phi]d\sigma\\
\nm\ds &-\frac{\ep}{4\pi} \int_{\p \Om}\Big[-2\Scal_{\Om}[\tau h\phi]+\Dcal_{\Om}[ h\phi] \big|_{+}\Big]\pd{\Scal_{\Om}[\phi]}{n}\Big|_{+}d\sigma-\frac{\ep}{4\pi} \int_{\p \Om} h\Big[\pd{\Scal_{\Om}[\phi]}{n}\Big|_{+} \Big]^2 d\sigma\\
\nm \ds &+\frac{\ep}{4\pi} \int_{\p \Om}\frac{1}{\sqrt{det(\mathcal{G}})}\left(\nabla_{\xi,\ta} \cdot
\Big(h\sqrt{det(\mathcal{G})}\mathcal{G}^{-1}\nabla_{\xi,\ta}\Scal_{\Om}[\phi]\Big)\right)\Scal_{\Om}[\phi]d\sigma+O(\ep^2),
\end{align*}
where $\phi:=\phi_\ep\circ \Psi_\ep$.  From the decomposition of $\phi$ in \eqref{decomposition-phi-ep}, we write
\begin{align*}
\ds \textrm{cap}(\Om_\ep)=&-\frac{1}{4\pi} \int_{\p \Om}\pd{\Scal_{\Om}[\phi_0]}{n}\Big|_{+} \Scal_{\Om}[\phi_0]d\sigma-\frac{\ep}{4\pi} \int_{\p \Om} h\Big[\pd{\Scal_{\Om}[\phi_0]}{n}\Big|_{+} \Big]^2 d\sigma\\
\nm \ds &-\frac{\ep}{4\pi} \int_{\p \Om}\bigg[\pd{\Scal_{\Om}[\phi_1]}{n}\Big|_{+}-2\pd{\Scal_{\Om}[\tau h\phi_0]}{n}\Big|_{+}+\pd{\Dcal_{\Om}[ h\phi_0]}{n}  \bigg]\Scal_{\Om}[\phi_0]d\sigma\\
\nm\ds &-\frac{\ep}{4\pi} \int_{\p \Om}\Big[\Scal_{\Om}[\phi_1]-2\Scal_{\Om}[\tau h\phi_0]+\Dcal_{\Om}[ h\phi_0]\big|_{+} \Big]\pd{\Scal_{\Om}[\phi_0]}{n}\Big|_{+}d\sigma\\
\nm \ds &+\frac{\ep}{4\pi} \int_{\p \Om}\frac{1}{\sqrt{det(\mathcal{G}})}\left(\nabla_{\xi,\ta} \cdot
\Big(h\sqrt{det(\mathcal{G})}\mathcal{G}^{-1}\nabla_{\xi,\ta}\Scal_{\Om}[\phi_0]\Big)\right)\Scal_{\Om}[\phi_0]d\sigma+O(\ep^2).
\end{align*}
Since  $\Scal_{\Om}[\phi_0]=1$  on $\p \Om$, we get  $\nabla_{\xi,\ta}\Scal_{\Om}[\phi_0]=0$  on $\p \Om$  and then the last  integral is equal to zero.   By using Green's formula, we deduce that the third  integral is equal to the forth integral.
Therefore
\begin{align}\label{5}
\ds \textrm{cap}(\Om_\ep)=&\textrm{cap}(\Om)-\frac{\ep}{4\pi} \int_{\p \Om} h\Big[\pd{\Scal_{\Om}[\phi_0]}{n}\Big|_{+} \Big]^2 d\sigma\nonumber\\
\nm\ds &-\frac{\ep}{2\pi} \int_{\p \Om}\Big[\Scal_{\Om}[\phi_1]-2\Scal_{\Om}[\tau h\phi_0]+\Dcal_{\Om}[ h\phi_0] \big|_{+} \Big]\pd{\Scal_{\Om}[\phi_0]}{n}\Big|_{+}d\sigma+O(\ep^2).
\end{align}
It follows from \eqref{444} that
\begin{align}\label{6}
\ds \int_{\p \Om}\Big[\Scal_{\Om}[\phi_1]-2\Scal_{\Om}[\tau h\phi_0]+\Dcal_{\Om}[ h\phi_0] \big|_{+} \Big]\pd{\Scal_{\Om}[\phi_0]}{n}\Big|_{+}d\sigma=-\int_{\p \Om} h\Big[\pd{\Scal_{\Om}[\phi_0]}{n}\Big|_{+} \Big]^2 d\sigma+O(\ep).
\end{align}
Finally, we conclude from the representation formula of $u$ \eqref{representation-u}, \eqref{5}, and \eqref{6} that
\begin{align*}
\ds \textrm{cap}(\Om_\ep)=&\textrm{cap}(\Om)+\frac{\ep}{4\pi} \int_{\p \Om} h\Big[\pd{u}{n} \Big]^2 d\sigma+O(\ep^2).
\end{align*}
This completes the proof of Theorem \ref{main-Theorem}, as desired.

\subsubsection{Proof of Theorem \ref{main-theorem-3}}\label{proof-theorem-3}
By \eqref{representation-u-ep} and \eqref{representation-u}, we have
\begin{align*}
\ds u_\ep(x)-u(x)=-\frac{1}{4\pi}\int_{\p \Om} \Gamma(x-\tilde y) \pd{u_{\ep}}{\tilde{n}}(\tilde y)d\tilde\sigma(\tilde y)
+\frac{1}{4\pi}\int_{\p \Om} \Gamma(x-y) \pd{u}{{n}}(y)d\sigma( y).
\end{align*}
It then follows from Theorem \ref{main-Theorem-current}, \eqref{dsigma-ep}, and \eqref{asymptotic-green-function-infinity} that
\begin{align*}
\ds u_\ep(x)-u(x)=&-\frac{\ep}{4\pi}\int_{\p \Om} \Gamma(x-\tilde y) \pd{v}{{n}}(y)d\sigma( y)
\\
\nm \ds  & -\frac{1}{4\pi}\int_{\p \Om} \Big[\Gamma(x-\tilde y)-\Gamma(x-y) \Big]\pd{u}{{n}}(y) d\sigma( y)+O(\frac{\ep^2}{|x|}).
\end{align*}
Since
\begin{align*}
\ds \Gamma(x-\tilde y)=\frac{1}{|x|}+ O(\frac{1}{|x|^2}), \q\q \Gamma(x-\tilde y)-\Gamma(x-y) = O(\frac{\ep}{|x|^2}) \q \mbox{as } |x|\rightarrow +\infty.
\end{align*}
Therefore
\begin{align*}
\ds u_\ep(x)-u(x)= -\frac{\ep}{4\pi|x|}\int_{\p \Om} \pd{v}{{n}}(y)d\sigma( y)+O(\frac{\ep}{|x|^2})+O(\frac{\ep^2}{|x|}).
\end{align*}
According to  the Green's formula, we immediately see that
\begin{align*}
\ds \int_{\p \Om} \pd{v}{{n}}d\sigma=\int_{\p \Om} \pd{v}{{n}}u d\sigma=\int_{\p \Om} v \pd{u}{{n}} d\sigma=-\int_{\p \Om} h \Big[\pd{u}{{n}}\Big]^2 d\sigma.
\end{align*}
This completes the proof of Theorem \ref{main-theorem-3}.

\section{Sensitivity analysis of the first eigenvector  of the   operator $\Kcal^{*}_{\Om}$  with respect to
 small perturbations in the surface  of   its domain  }

Let $\Om$ be a bounded domain in $\RR^3$ with $\mathcal{C}^2$-boundary $\p \Om$.  The spectrum of $\Kcal_{ \Om}^*: L^2(\p\Om)\rightarrow L^2(\p\Om)$  is discrete, lies in the interval $(-\frac{1}{2}, \frac{1}{2}]$,
and accumulates at zero.  More precisely, let   $\{\lambda_j\}_{0}^{\infty}$   be the  eigenvalues of $\Kcal_{ \Om}^*$ on $L^2(\p\Om)$,  then, the first eigenvalue $\lambda_0$ is equal to ${1}/{2}$ and has   geometric multiplicity $1$ while  $\lambda_j \in (-\frac{1}{2}, \frac{1}{2})$ for $j\geq 1$ with  $|\lambda_1|\geq |\lambda_2|\geq \cdots \rightarrow 0$ as $j \rightarrow \infty$ arranged  repeatedly according to their multiplicities. See for example \cite{book-3}.

Denote by  $\varphi_0$  the  first eigenvector of $\Kcal_{ \Om}^*$  on $L^2(\p \Om)$ associated to the first  eigenvalue ${1}/{2}$ with $\|\varphi_0\|_{L^2(\p \Om)}=1$.
We claim that  $\varphi_0$ is equal to $\frac{\p u}{\p n}{/}\|\frac{\p u}{\p n}\|_{L^2(\p \Om)}$, where  $u$    represents the electrostatic potential in the presence of  the conductor $\Om$ in  electrostatic equilibrium, it is
the unique solution of \eqref{u}. In fact, it follows from \eqref{representation-u}  and \eqref{nuS} that
\begin{align*}
\ds \pd{u}{n}=\Big(\frac{1}{2}I+\Kcal_{\Om}^*\Big)\Big[\pd{u}{n}\Big]\q \mbox{on } \p\Om,
\end{align*}
namely,
\begin{align}
\ds \Kcal_{\Om}^*\Big[\pd{u}{n}\Big]=\frac{1}{2}\pd{u}{n}\q \mbox{on } \p\Om.\label{relationship-varphi-u}
\end{align}
From the uniqueness of   $\varphi_0$, we deduce that $ \varphi_0=\frac{\p u}{\p n}{/}\|\frac{\p u}{\p n}\|_{L^2(\p \Om)}$ on  $\p\Om.
$

It is known that the first eigenvalue ${1}/{2}$
is independent of $\p \Om$, that is, it does not  affected by any   smooth perturbations of $\p \Om$. In view of this  remark,
the electrostatic capacity of an isolated  conductor may also  be defined as the amount of a charge required to raise
 the first eigenvalue of the $\Kcal^{*}_{\Om}$  operator   at ${1}/{2}$.

 Similarly to \eqref{relationship-varphi-u}, we have
 \begin{align*}
\ds \Kcal_{\Om_\ep}^*\Big[\pd{u_\ep}{\tilde{n}}\Big]=\frac{1}{2}\pd{u_\ep}{\tilde{n}}\q \mbox{on } \p\Om_\ep,
\end{align*}
 where  $u_{\ep}$  is the unique  solution
 of  \eqref{uep}. Therefore  the first eigenvector $\varphi^\ep_0$  of the operator  $\Kcal_{\p \Om_\ep}^*$ on $L^2(\p \Om_\ep)$  with the eigenvalue $1/2$
is equal to $\frac{\p u_\ep}{\p \tilde{n}}{/}\|\frac{\p u_\ep}{\p \tilde{n}}\|_{L^2(\p \Om_\ep)} $.

From Theorem \ref{main-Theorem-current}, we obtain  in  the following theorem the fourth result of this paper,  an  asymptotic expansion for the first eigenvector $\varphi_0^\ep$ on
 $\p \Om_\ep$    as $\ep \rightarrow 0$.
\begin{thm} \label{main-Theorem-eigenvector} Let $\tilde {x}=x+\ep h(x)n(x)\in \p \Om_{\ep}$, for $x\in \p \Om$.  Let $\varphi_0^\ep$  and $\varphi_0$ be the first eigenvectors of $\Kcal_{\Om_\ep}^*$  and $\Kcal_{ \Om}^*$ with the eigenvalue ${1}/{2}$, respectively.
 The following  asymptotic  expansion  holds:
\begin{align*}
\ds   \varphi_0^\ep(\tilde x)=&\varphi_0(x)+2\ep\tau(x) h(x) \varphi_0(x)+\ep \tilde{v}(x)-\ep \la \tau h \varphi_0+\tilde{v}, \varphi_0 \ra\varphi_0(x)+O(\ep^{2}),
\end{align*}\label{asymptotic-eigenvector}
with $ \tilde{v}=\frac{\p v}{\p n} /\|\frac{\p u}{\p n}\|_{L^2(\p \Om)}$, where  $u$ and $v$  are the unique solutions  of \eqref{u} and \eqref{v}, respectively,
and the remainder $O(\ep^{2})$ depends  only
on the $\mathcal{C}^2$-norm of $X$ and the $\mathcal{C}^1$-norm of $h$.
\end{thm}

The  fifth  result of this paper is the following theorem, an asymptotic expansion of $ \int_{\p \Om} \big(\varphi_0^\ep(\tilde x) - \varphi_0 (x) \big )\varphi_0 (x)d{\sigma}(x)$  as $\ep \rightarrow 0$.
\begin{thm} \label{main-Theorem-integral-eigenvector} Let $\varphi_0^\ep$  and $\varphi_0$ be the first eigenvectors of $\Kcal_{ \Om_\ep}^*$  and $\Kcal_{\Om}^*$ with the eigenvalue ${1}/{2}$, respectively. The following  asymptotic  expansion  holds:
\begin{align}
\ds  \int_{\p \Om} \big(\varphi_0^\ep(\tilde x) - \varphi_0 (x) \big )\varphi_0 (x)d{\sigma}(x)=& \ep  \int_{\p \Omega}\tau(x) h(x) [\varphi_0(x)]^2 d\sigma(x) +O(\ep^2),
\label{main-Theorem-equality-eigenvector}
\end{align}
where   the remainder $O(\ep^2)$ depends  only
on the $\mathcal{C}^2$-norm of $X$ and the $\mathcal{C}^1$-norm of $h$.
\end{thm}


The asymptotic expansion \eqref{main-Theorem-equality-eigenvector} could be used to determine some
properties  on the shape perturbation of an object   from measurements on the perturbed shape itself (see \cite {Z2}) of the first eigenvector of the $L^2$-adjoint  of the NP   operator.



\begin{thebibliography}{10}

\bibitem{Book-Ammari}  H. Ammari, \textit{An Introduction to Mathematics of Emerging Biomedical Imaging}, Math.
Appl., Volume 62, Springer, Berlin, 2008.

\bibitem{ABFKL} H. Ammari, E. Beretta, E. Francini, H. Kang, and M. Lim,
Optimization algorithm for reconstruction interface changes of a conductivity inclusion
from modal measurements, Math. comp., 79 (2010), 1757-1777.


\bibitem{book-3}  H. Ammari, B. Fitzpatrick, H. Kang, M. Ruiz, S. Yu, and H. Zhang,\textit{ Mathematical and Computational Methods in Photonics and Phononics},  Mathematical Surveys and Monographs, Volume 235, American Mathematical Society, Providence, 2018.


\bibitem{book-1} H. Ammari and H. Kang \textit{Polarization and Moment Tensors with Applications to Inverse Problems
and Effective Medium Theory}, Applied Mathematical Sciences, Vol. 162, SpringerVerlag, New York, 2007.


\bibitem{AKLZ1} H. Ammari, H. Kang, M. Lim, and H. Zribi,
Conductivity interface problems. Part I: small perturbations of an
interface,  Trans. Amer. Math. Soc., 362 (2010), 2435-2449.

\bibitem{AKLZ2} H. Ammari, H. Kang, M. Lim, and H. Zribi,
The generalized polarization tensors for resolved imaging. Part I:
shape reconstruction of a conductiovity inclusion,  Math. of  comp., 81 (2012), 367-386.



\bibitem{CMM82}
{R.R. Coifman, A. McIntosh, and Y. Meyer}, { L'int{\'e}grale de
Cauchy d{\'e}finit un op{\'e}rateur born{\'e} sur $L^2$ pour les
courbes lipschitziennes}, Ann. Math., 116 (1982), 361--387.

\bibitem{Coifman} R. Coifman, M. Goldberg, T. Hrycak, M. Israeli, and V. Rokhlin,
 \newblock An improved operator expantion algorithm for direct
and inverse scattering computations,
 \newblock  Waves Random Media, 9 (1999), 441-457.


 \bibitem{Folland76} G.B. Folland, \textsl{Introduction to Partial Differential
Equations}, Princeton University Press, Princeton, New Jersey, 1976.

 \bibitem{Jerison} D. Jerison, A Minkowski problem for electrostatic capacity, Acta Math., 176 (1996) 1-47.

\bibitem{Kellogg} O.D. Kellogg, \textit{ Foundations of Potential Theory}, Dover, New York, 1953.

\bibitem{KZ1} A. Khelifi and H. Zribi, Asymptotic expansions for the voltage potentials with two- and three-dimensional thin interfaces,
Math.  Methods. Appl. Sci., 34 (2011), 2274-2290.

 \bibitem{KZ2} A. Khelifi and  H. Zribi, Boundary voltage perturbations resulting from small
surface changes of a conductivity  inclusion, Appl. Anal., Vol.  93 (2014), 46-64.

\bibitem{LTZ} J. Lagha,  F. Triki, and H. Zribi, Small perturbations of an interface for elastostatic problems,
Math.  Methods. Appl. Sci.,  40 (10)(2017), 3608-3636.

\bibitem{LZ} J. Lagha  and H. Zribi, An asymptotic expansion for perturbations in the displacement field due to the presence of
thin interfaces, Appl. Anal., volume 1, (2017), 1-23.


\bibitem{LLZ} M. Lim, K. Louati and H. Zribi,
\newblock Reconstructing small perturbations of scatterers from electric or acoustic far‐field measurements,
\newblock  Math.  Methods. Appl. Sci., 31 (2008), no 11,  1315-1332.

\bibitem{volume2} S. J. Ling, J. Sanny, and W. Moebs, \textit{University Physics - Volume 2} (OpenStax), (2016).


\bibitem{Maxwell}  J. C. Maxwell, \textit{An elementary treatise on electricity},  Clarendon Press in Oxford, 1881.



\bibitem{Poincare} H. Poincar\'e, Figures d'\'equilibre d'une masse fluide, Paris, 1902.

  \bibitem{Polya}  G. Polya, Estimating electrostatic capacity,  Am. Math. Mon.,  54, no. 4 (1947),  201–206.

  \bibitem{Polya-Szego-book}  G. P\'olya and G. Szeg\"o, \textit{Isoperimetric Inequalities in Mathematical Physics, Annals
of Mathematical Studies}, Number 27, Princeton University Press, Princeton, NJ, 1951.



\bibitem{Szego} G. Szeg\"o,  On the capacity of a condenser, Bull. Amer. Math. Soc. vol. 51 (1945) pp. 325-350.



 \bibitem{Z}  H. Zribi, Asymptotic expansions for currents caused by small interface changes of an electromagnetic inclusion, Appl. Anal., 92, (2013), 172-190.

 \bibitem{Z2}  H. Zribi,  Reconstructing  small perturbations of an  obstacle
 for acoustic waves from boundary measurements  on the perturbed shape itself,  Math.  Methods. Appl. Sci.,  45 (2022), no 1,  93-112.



















\end{thebibliography}
\end{document}